\documentclass[a4paper]{article}
\usepackage{amsmath}
\usepackage{amssymb}
\usepackage{amsfonts}
\usepackage{theorem}
\usepackage{blkarray}
\usepackage{arydshln}
\usepackage{multirow}
\usepackage{graphicx}

\makeatletter
\@addtoreset{equation}{section}

\makeatother
\newtheorem{theo}{Theorem}[section]
\newtheorem{defi}{Definition}[section]
\newtheorem{lemm}{Lemma}[section]
\newtheorem{prop}{Proposition}[section]

\newtheorem{rema}{Remark}[section]

\newtheorem{claim}{Claim}[section]
\newcommand{\Z}{\mathbb Z}
\newcommand{\C}{\mathbb C}

\newcommand{\R}{\mathbb R}
\newcommand{\bP}{\mathbb P}
\makeatletter
\newcommand{\xRightarrow}[2][]{%
\ext@arrow 0055{\Rightarrowfill@}{#1}{#2}%
}
\def\Rightarrowfill@{\arrowfill@\Relbar\Relbar\Rightarrow}
\newcommand{\xLeftarrow}[2][]{%
\ext@arrow 0055{\Leftarrowfill@}{#1}{#2}%
}
\def\Leftarrowfill@{\arrowfill@\Leftarrow\Relbar\Relbar}
\newcommand{\xLongleftrightarrow}[2][]{
\ext@arrow 0055{\llrafill@}{#1}{#2}}
\def\llrafill@{\arrowfill@\Leftarrow\Relbar\Rightarrow}

\begin{document}
\title{
\begin{flushright}
  \begin{minipage}[b]{5em}
    \normalsize
    ${}$      \\
  \end{minipage}
\end{flushright}
{\bf Moduli Space of Quasi-Maps from $\mathbb P^{1}$  with Two Marked Points to $\mathbb P(1,1,1,3)$ and $j$-invariant}}
\author{Masao Jinzenji${}^{1)}$, Hayato Saito${}^{2)}$\\
\\
\it Division of Mathematics, Graduate School of Science \\
\it Hokkaido University \\
\it  Kita-ku, Sapporo, 060-0810, Japan\\
{\it e-mail address: 1) jin@math.sci.hokudai.ac.jp}\\
\hspace{3.5cm}{\it                           2) hayato@math.sci.hokudai.ac.jp }}
\maketitle

\begin{abstract}
In this paper, we construct toric data of moduli space of quasi maps of degree $d$ from $\bP^{1}$ with two marked points  
to weighted projective space $\bP(1.1,1,3)$. With this result, we prove that the moduli space is a compact toric orbifold.  
We also determine its Chow ring. Moreover, we give a proof of the conjecture proposed by Jinzenji that a series of intersection
numbers of the moduli spaces coincides with expansion coefficients of inverse function of $-\log(j(\tau))$.  
\end{abstract}

\section{Introduction.}
\subsection{Overview of Classical Mirror Symmetry}
\textit{Mirror Symmetry} is a symmetry between two topological sigma models, A-model on a Calabi-Yau manifold $X$ and 
B-model on another Calabi-Yau manifold $X^{*}$.
Mathematically, A-model has information of world sheet instantons, i.e., holomorphic maps from a Riemann surface $\Sigma$
to $X$.
On the other hand, B-model has information of deformation of Hodge structure of $X^{*}$.
Mirror Symmetry surprisingly connects these concepts and has fascinated both physicists and mathematicians.

Classical Mirror Symmetry enables us to compute GW invariants of $X$ by using solutions of Picard-Fuchs equations for
Period integral of $X^{*}$ (\cite{CdGP, CK, Givental}).
GW invariants are defined as intersection numbers of the moduli spaces of stable maps $\overline{M}_{g,n}(X,\beta)$.
But computation of the invariants from this definition are usually done by localization technique and the process 
is quite complicated. 
On the other hand, the Picard-Fuchs equations of $X^{*}$ are linear differential equations and easy to solve in many cases.
Furthermore, \textit{mirror map} which connects A-model and B-model is also given by solutions of Picard-Fuchs equations.
Hence, the process of computation of GW invariants by using classical mirror symmetry is much simpler than the process by
using localization.

In \cite{Jin1}, one of the authors (Jinzenji) introduced moduli space $\widetilde{Mp}_{0,2}(N,d)$, which 
is a compactified moduli space of quasi maps (polynomial maps) from $CP^{1}$ to $CP^{N-1}$ with two marked points.
This moduli space is deeply connected to informations of the B-model in classical mirror symmetry.
Namely, generating functions of intersection numbers on $\widetilde{Mp}_{0,2}(N,d)$ directly gives us the solution of Picard-Fuchs equations and the mirror map.
The geometrical feature of $\widetilde{Mp}_{0,2}(N,d)$ comes from the fact that its geometrical structure is much simpler 
than the corresponding moduli space $\overline{M}_{0,2}(CP^{N-1},d)$, which is compactified by using stable maps.  
In \cite{S}, another author (Saito) constructed a concrete toric data of $\widetilde{Mp}_{0,2}(N,d)$.
This toric data gives us a lot of properties of $\widetilde{Mp}_{0,2}(N,d)$.
For example, $\widetilde{Mp}_{0,2}(N,d)$ is a compact toric orbifold.
We can also determine Chow ring of $\widetilde{Mp}_{0,2}(N,d)$ from the toric data.
In \cite{S}, Saito also discovered an injective homomorphism from $A^*(\widetilde{Mp}_{0,2}(N,d))$ to $A^*(\overline{M}_{0,2}(\bP^{N-1},d))$ for $d=1,2$ cases.
Furthermore, Saito proved a formula that describes genus $0$ GW invariants of projective hypersurfaces of $d=1,2$ in terms of  
Chow ring of $\overline{M}_{0,2}(\bP^{N-1},d)$.
This formula was found by pursuing analogy of the formula that describes the corresponding intersection numbers of $\widetilde{Mp}_{0,2}(N,d)$ in terms of its Chow ring, that was implied in \cite{Jin1}.

\subsection{$j$-invariant.}
The $j$-invariant is a weight zero modular function of $\tau$ which is the coordinate of complex structures of elliptic curves over $\C$:
\begin{equation*}
 \C/(\Z\oplus\Z\tau),
\end{equation*}
where $\tau$ is in the upper half plane.
$j$-invariant is one of the fundamental tools for studying elliptic curves.
For example, it determines a group structure of elliptic curves.

Let us introduce $q=e^{2\pi \sqrt{-1}\tau}$, then Fourier expansion of j-invariant is given by,
\begin{align}
 j(\tau)&=q^{-1}+744+196884q+21493760q^2+864299970q^3+20245856256q^4+\cdots\nonumber\\
 &=:q^{-1}+\sum_{d=1}^\infty j_d\cdot q^{d-1}.
\end{align}
These coefficients were found to be related to ranks of irreducible representations of the Monster group (the largest sporadic simple group), which is known as monstrous moonshine.

In this paper, we deal with the expansion coefficients of \textit{inverse} function of $\displaystyle{-\log(j(\tau))}$,
\begin{align}
 2\pi\sqrt{-1}\tau&= -\log(j)+744j^{-1}+473652j^{-2}+451734080j^{-3}+510531007770j^{-4}+\cdots\nonumber\\
 &=:-\log(j)+\sum_{d=1}^\infty w_d \cdot j^{-d}.
\end{align}
This function appeared as the mirror map of the K3 surface in $\mathbb{P}(1,1,1,3)$. 
The expansion coefficient $j_{d}$ is reconstructed by the expansion coefficient $w_{d}$ via the following
relation:
\begin{equation}
j_{d}=\sum_{\sigma_d\in OP_d}(-(d-1))^{l(\sigma_d)-1}\frac{1}{(l(\sigma_d))!}
\prod_{j=1}^{l(\sigma_{d})}w_{d_j},
\end{equation}
where $OP_{d}$ is set of ordered partitions of positive integer $d$,
\begin{equation*}
OP_{d}:=\{\sigma_{d}=(d_{1},d_{2},\cdots,d_{l(\sigma_{d})})\;|\;\sum_{j=1}^{l(\sigma_{d})}d_{j}=d,\;\;d_{1},d_{2},\cdots,d_{l(\sigma_{d})}\geq 1\;\},
\end{equation*}
and $l(\sigma_{d})$ is length of an ordered partition $\sigma_{d}\in OP_{d}$.
In \cite{Jin1} (arXiv version), Jinzenji conjectured that {\bf the coefficient $w_{d}$ is written as an intersection number 
of moduli space of quasi maps from $\mathbb{P}^{1}$ to $\mathbb{P}(1,1,1,3)$ of degree $d$ with two marked points}.   
The aim of our paper is to prove this conjecture. 

\subsection{Picard-Fuchs equation for $j$-invariant.}\label{subsec:PF eq for j-inv}
It is known that the Picard-Fuchs equation for 1-parameter deformation of algebraic K3 surfaces is solved in terms of the j-invariant (\cite{LY}).
Let us demonstrate it in the case of 1-parameter deformation of algebraic K3 surface in $\bP(1,1,1,3)$.
We mainly refer section 5.4. of \cite{LY}.

The 1-parameter deformation of K3surface which we use is given by,
\begin{equation*}
 x_1^6+x_2^6+x_3^6+x_4^2+z^{-1/6}x_1x_2x_3x_4=0,
\end{equation*}
embedded in  $\bP(1,1,1,3)=\{(x1,x2,x3,x4)\}$.
Its Picard-Fuchs equation is given by,
\begin{equation}\label{eq:P-F eq of K3}
 (\Theta^3-8z(6\Theta+1)(6\Theta+3)(6\Theta+5))f(z)=0,
\end{equation}
where $\Theta$ is a differential operator $\Theta=z\frac{d}{dz}$.
It is obtained from standard techniques, for example, Griffiths-Dwork method, or {\cal A}-Hypergeometric equations (see \cite{CK}, etc.).

In the following, we briefly review the process of solving the above equation. 
 Let us assume that the solution $f(z)$ is expanded as follows:
\begin{equation*}
 f(z)=\sum_{n=0}^\infty a_n(\epsilon)z^{n+\epsilon},
\end{equation*}
where $\epsilon$ is a parameter and $a_0(\epsilon)=1$..
By substituting it to (\ref{eq:P-F eq of K3}), we obtain,
\begin{equation*}
 \sum_{n=0}^\infty (a_n(\epsilon)(n+\epsilon)^3-8(6n-5+6\epsilon)(6n-3+6\epsilon)(6n-1+6\epsilon)a_{n-1}(\epsilon))z^{n+\epsilon}=0,
\end{equation*}
where $a_{-1}(\epsilon)=0$.
Since it holds for any $z$, we obtain, 
\begin{align*}
 a_n(\epsilon)&=\frac{8(6n-5+6\epsilon)(6n-3+6\epsilon)(6n-1+6\epsilon)}{(n+\epsilon)^3}a_{n-1}(\epsilon)\\
 &=\frac{\Gamma(6n+6\epsilon+1)}{\Gamma(n+\epsilon+1)^3\Gamma(3n+3\epsilon+1)},
\end{align*}
where $\Gamma(x)$ is the Gamma function.
Therfore, we can obtain the solutions of (\ref{eq:P-F eq of K3}) from the following expression:
\begin{equation*}
 f(z,\epsilon):=\sum_{n=0}^\infty \frac{\Gamma(6n+6\epsilon+1)}{\Gamma(n+\epsilon+1)^3\Gamma(3n+3\epsilon+1)}z^{n+\epsilon}.
\end{equation*}
.
By setting $\epsilon=0$, it gives a solution which is holomorphic at $z=0$:
\begin{equation}\label{eq:solution0}
 f_0(z):=\sum_{d=0}^\infty \frac{2^{3d}(6d-1)!!}{(d!)^3}z^d.
\end{equation}
If we differentiate $f(z,\epsilon)$ by $\epsilon$ and set $\epsilon=0$, we obtain the solution which has a log-singularity at $z=0$:
\begin{equation}\label{eq:solution1}
 f_1(z):=f_0(z)(\log(z))+\sum_{d=0}^\infty (\sum_{j=1}^{3d}\frac{6}{2j-1}-\sum_{j=1}^{d}\frac{3}{j})\frac{2^{3d}(6d-1)!!}{(d!)^3}z^d.
\end{equation}

The mirror map for the K3 surface is given by,
\begin{equation*}
 f_1(z)/f_0(z).
\end{equation*}
It gives us the inverse function of $\displaystyle{-\log(j(\tau))}$ (\cite{LY}).
\begin{theo}[B. Lian, S. T. Yau, 1996]\label{theorem:LY}
 \begin{align*}
  \frac{f_1(j^{-1})}{f_0(j^{-1})}=2\pi\sqrt{-1}\tau=-\log(j)+\sum_{d=1}^\infty w_d\cdot j^{-d}.
 \end{align*}
\end{theo}

\subsection{The Goal of this paper.}
In order to prove the conjecture, we explicitly construct the compactified moduli space of degree $d$ quasi maps
from $\bP^{1}$ to $\bP(1,1,1,3)$ with two marked points, which we denote by $\widetilde{Mp}_{0,2}(\bP(1,1,1,3),d)$. 
In section \ref{section:P1113}, we provide a toric data of $\widetilde{Mp}_{0,2}(\bP(1,1,1,3),d)$ and prove 
the following theorem.
\begin{theo}\label{theo:properties MpP1113}
 $\widetilde{Mp}_{0,2}(\bP(1,1,1,3),d)$ is a compact toric orbifold.
\end{theo}
Furthermore, we show that the Chow ring of $\widetilde{Mp}_{0,2}(\bP(1,1,1,3),d)$ is given by,
\begin{prop}\label{prop:Chow ring of MpP1113}
 $A^*(\widetilde{Mp}_{0,2}(\bP(1,1,1,3),d))\cong \C [H_0,H_1,\dots,H_d]/{\cal I}$,\\
where ${\cal I}=(H_0^4(2H_0+H_1),H_1^4(H_0+2H_1)(2H_1+H_2)(-H_0+2H_1-H_2),\dots,H_{d-1}^4(H_{d-2}+2H_{d-1})(2H_{d-1}+H_d)(-H_{d-2}+2H_{d-1}-H_d),H_d^4(H_{d-1}+2H_d))$.
\end{prop}
With these preparations, we define the intersection number that corresponds to $w_{d}$ by,
\begin{defi}\label{definition:intersection num of MpP1113}
 \begin{equation*}
  w({\cal O}_{z^a}{\cal O}_{z^b})_{0,d}:=\int_{\widetilde{Mp}_{0,2}(\bP(1,1,1,3),d)} H_0^aH_1^b \cdot \frac{\displaystyle{\prod_{i=1}^de^6(H_{i-1},H_i)}}{\displaystyle{\prod_{i=1}^{d-1}6H_i}},
 \end{equation*}
 where
 \begin{equation*}
  e^6(x,y):=\prod_{j=0}^6 ((6-j)x+y).
 \end{equation*}
\end{defi}
In this definition, $H_0,H_1,\dots,H_d$ are generators of Chow rings of $\widetilde{Mp}_{0,2}(\bP(1,1,1,3),d)$.

In section \ref{section:j-inv}, we prove our main result.
\begin{theo}[Main Theorem]\label{theorem:main}
 \begin{equation}\label{eq:j-inv conj}
  w_d=\frac12w({\cal O}_{z^1}{\cal O}_{z^0})_{0,d}.
 \end{equation}
\end{theo} 

Of course, our main result also follows from general theory of wall crossing formula of quasimap theory by Cheong, Ciocan-Fontanine 
and Kim \cite{CK, CCK}. But it seems that their theory does not supply explicit information of Chow ring of our moduli space of 
quasimaps from $\mathbb P^{1}$  with two marked points to $\mathbb P(1,1,1,3)$.  Our approach in this paper is quite direct 
and provides explicit computation of intersection numbers that correspond to expansion coefficients of the mirror map. 
Therefore, we expact that our approach also contributes to deeper understanding of machinery of mirror computation.

\section{The moduli space $\widetilde{Mp}_{0,2}(\bP(1,1,1,3),d)$.}\label{section:P1113}
In this section, we provide definition of the quasi-map moduli space $\widetilde{Mp}_{0,2}(\bP(1,1,1,3),d)$ and prove the
following theorem.
\begin{theo}
 $\widetilde{Mp}_{0,2}(\bP(1,1,1,3),d)$ is a compact toric orbifold.
\end{theo}
Furthermore, we compute its Chow ring:
\begin{equation*}
  A^*(\widetilde{Mp}_{0,2}(\bP(1,1,1,3),d))\cong \C [H_0,H_1,\dots,H_d]/{\cal I},
\end{equation*}
where ${\cal I}=(H_0^4(2H_0+H_1),H_{1}^4(H_{0}+2H_{1})(2H_{1}+H_{2})(-H_{0}+2H_{1}-H_{2}),
\cdots,H_{j}^4(H_{j-1}+2H_{j})(2H_{j}+H_{j+1})(-H_{j-1}+2H_{j}-H_{j+1}),\cdots,
H_{d-1}^4(H_{d-2}+2H_{d-1})(2H_{d-1}+H_{d})(-H_{d-2}+2H_{d-1}-H_{d}),
H_d^4(H_{d-1}+2H_d))$.

In addition, we define an intersection numbers $w({\cal O}_{z^a}{\cal O}_{z^b})_{0,d}$ of the quasi-map moduli space $\widetilde{Mp}_{0,2}(\bP(1,1,1,3),d)$.

\subsection{The fan of $\widetilde{Mp}_{0,2}(\bP(1,1,1,3),d)$.}
The generic quasi-map from $\bP^1$ to $\bP(1,1,1,3)$ is given by
\begin{align*}
 &\bP^1\rightarrow \bP(1,1,1,3) \\
 &[s:t]\mapsto [f_0(s.t):f_1(s,t):f_2(s,t):f_3(s,t)],
\end{align*}
where
\begin{align*}
 &f_i(s,t):=\sum_{j=0}^d a_{i,j}s^{d-j}t^j,& (0\leq i\leq 2),\\
 &f_3(s,t):=\sum_{j=0}^{3d} a_{3,j} s^{3d-j}t^j. &
\end{align*}
The following $(\C^*)^2$-action on $(a_{i,j})$ is induced from projective equivalence of $\bP(1,1,1,3)$ and 
automorphism group of $\bP^{1}$ which keeps $(0:1), (1:0)\in \bP^1$ fixed.
\begin{align}\label{eq:action pP1113}
 &(\mu,\nu)\cdot ({\bf a}_0,{\bf a}_1,\dots,{\bf a}_d,a_{3,0},a_{3,1},\dots,a_{3,3d})\nonumber\\
 =&(\mu {\bf a}_0, \mu\nu {\bf a}_1,\mu\nu^2 {\bf a}_2,\dots,\mu\nu^d{\bf a}_d, \mu^3a_{3,0},\mu^3\nu a_{3,1},\mu^3\nu^2a_{3,2},\dots,\mu^3\nu^{3d}a_{3,3d}).
\end{align}
Here ${\bf a}_{j}\in {\bf C}^3$ represents a vector $(a_{0,j},a_{1,j},a_{2,j})$.
This is equivalent to the following action:
\begin{align}\label{eq:action P1113}
 &(\mu',\nu')\cdot ({\bf a}_0,{\bf a}_1,\dots,{\bf a}_d,a_{3,0},a_{3,1},\dots,a_{3,3d})\nonumber\\
 =&((\mu')^d {\bf a}_0, (\mu')^{d-1}(\nu') {\bf a}_1,(\mu')^{d-2}(\nu')^2 {\bf a}_2,\dots,(\nu')^d{\bf a}_d,\nonumber \\
 & (\mu')^{3d}a_{3,0},(\mu')^{3d-1}(\nu') a_{3,1},(\mu')^{3d-2}(\nu')^2a_{3,2},\dots,(\nu')^{3d}a_{3,3d}).
\end{align}

The set
\begin{align*}
 &Mp_{0,2}(\bP(1,1,1,3),d)\\
 :=&\{({\bf a}_0,{\bf a}_1,\dots,{\bf a}_d,a_{3,0},a_{3,1},\dots,a_{3,3d})\in \C^{3(d+1)+3d+1}\}/(\C^*)^2
\end{align*}
is not compact.
In order to see it, let us use the $(\C^*)^2$-action in (\ref{eq:action pP1113}) to turn $({\bf a}_0,a_{3,0})$ and $({\bf a}_d,a_{3,3d})$ into points in $\bP(1,1,1,3)$, $[{\bf a}_0,a_{3,0}]$ and $[{\bf a}_d,a_{3,3d}]$.
Then, we obtain
\begin{align*}
 &Mp_{0,2}(\bP(1,1,1,3),d) \\
 \cong& \{([{\bf a}_0,a_{3,0}],{\bf a}_1,{\bf a}_2,\dots,{\bf a}_{d-1},[{\bf a}_d,a_{3,3d}]\\
 &\quad ,a_{3,1},a_{3,2},\dots,a_{3,3d-1})\in\bP(1,1,1,3)\times\C^{3(d-1)+3d-1}\times\bP(1,1,1,3)\}/\Z_d.
\end{align*}
Here the $\Z_d$-action is given by
\begin{align*}
 ([{\bf a}_0,a_{3,0}],\zeta_d{\bf a}_1,\zeta_d^2{\bf a}_2,\dots,\zeta_d^{d-1}{\bf a}_{d-1},[{\bf a}_d,a_{3,3d}],\zeta_da_{3,1},\zeta_d^2a_{3,2},\dots,\zeta_d^{3d-1}a_{3,3d-1}),
\end{align*}
where $\zeta_d$ is the $d$-th primitive root of unity.

In order to compactify $Mp(\bP(1,1,1,3),d)$, we add variables $u_1,u_2,\dots,u_{d-1}$.
They should be added as in the case of $\widetilde{Mp}_{0,2}(N,d)$ (see \cite{Jin1} or \cite{S}).
Hence, we want to obtain a toric data of
\begin{align*}
 &\widetilde{Mp}_{0,2}(\bP(1,1,1,3),d)\\
 :=&\{ ({\bf a}_0,{\bf a}_1,\dots,{\bf a}_d,a_{3,0},a_{3,1},\dots,a_{3,3d},u_1,u_2,\dots,u_{d-1})\in U\}/(\C^*)^{d+1},
\end{align*}
where $U$ is a dense open subset of $\C^{3(d+1)+(3d+1)+(d-1)}$, and the $(\C^*)^{d+1}$-action is given by
\begin{align*}
 &(\lambda_0,\dots,\lambda_d)\cdot ({\bf a}_0,{\bf a}_1,\dots,{\bf a}_d,a_{3,0},a_{3,1},\dots,a_{3,3d},u_1,u_2,\dots,u_{d-1})\\
 =&(\lambda_0{\bf a}_0,\lambda_1{\bf a}_1,\dots,\lambda_d{\bf a}_d,\lambda_0^3a_{3,0},\lambda_0^2\lambda_1a_{3,1},\lambda_0\lambda_1^2a_{3,2},\lambda_1^3a_{3,3},\lambda_1^2\lambda_2a_{3,4},\dots,\lambda_d^3a_{3,3d},\\
&\lambda_0^{-1}\lambda_1^2\lambda_2^{-1}u_1,\lambda_1^{-1}\lambda_2^2\lambda_3^{-1}u_2,\dots,\lambda_{d-2}^{-1}\lambda_{d-1}^2\lambda_d^{-1}u_{d-1}).
\end{align*}
The structure of $U$ will be given by constructing toric data corresponding to $\widetilde{Mp}_{0,2}(\bP(1,1,1,3),d)$.

In the following, we construct a fan which is complete and simplicial, and realizes this $(\C^*)^{d+1}$-action.

Let $ p_0,p_1,p_2\in \Z^2$ be integer vectors given by,
\begin{equation*}
 (p_0,p_1,p_2)=\left(
 \begin{array}{cccccc}
  -1        &1         &0 \\
  -1        &0         &1 \\
 \end{array}
 \right).
\end{equation*}

Next, we introduce $(d+1)$ column vectors
\begin{align*}
 v'_0,v'_1,\dots,v'_d\in \Z^{d-1},
\end{align*}
defined by,
\begin{equation*}
  (v'_0,v'_1,\dots,v'_{d-1},v'_d)=\left(
 \begin{array}{ccccccc}
  -1       &2         &-1      &0         &\cdots  &0          &0 \\
  0        &-1        &2        &-1       & \cdots  &0         &0 \\
  0        &0         &-1       &2         &\cdots  &0         &0\\
  0        &0         &0        &-1        &\cdots  &0         &0\\
  \vdots &\vdots &\vdots &\vdots  &\ddots &\vdots   &0 \\
  0        &0         &0        &0         &\cdots  &-1        &0\\
  0        &0         &0        &0          &\cdots  &2         &-1
 \end{array}
 \right)\in M_{d-1,d+1}(\Z).
\end{equation*}
in the same way as the case of $\widetilde{Mp}_{0,2}(N,d)$ in \cite{S}.

In addtion, we have to introduce the following vectors:
\begin{align*}
 &(w_0,w_1,w_2,\dots,w_d):=\\
 &\left(
 \begin{array}{ccccccc}
  3       &0         &0        &0         &\cdots  &0        &0 \\
  2       &1         &0        &0         &\cdots  &0        &0 \\
  1       &2         &0        &0         &\cdots  &0        &0 \\
  0       &3         &0        &0         &\cdots  &0        &0 \\
  0       &2         &1        &0         &\cdots  &0        &0 \\
  0       &1         &2        &0         &\cdots  &0        &0 \\
  0       &0         &3        &0         &\cdots  &0        &0 \\
  0       &0         &2        &1         &\cdots  &0        &0 \\
  0       &0         &1        &2         &\cdots  &0        &0 \\
  0       &0         &0        &3         &\cdots  &0        &0 \\
  \vdots&\vdots&\vdots&\vdots&\ddots&\vdots&\vdots \\
  0       &0         &0        &0         &\cdots  &3        &0 \\
  0       &0         &0        &0         &\cdots  &2        &1 \\
  0       &0         &0        &0         &\cdots  &1        &2 \\
  0       &0         &0        &0         &\cdots  &0        &3 \\
 \end{array}
 \right)\in M_{3d+1,d+1}(\Z)
\end{align*}

Finally, we define column vectors,
\begin{align*}
 &v_{i,j}\;\;(0\leq i\leq 2, 0\leq j\leq d),\\
 &v_{3,j}\;\; (0\leq j\leq 3d),\\
 &u_k\;\; (1\leq k\leq d-1)
\end{align*}
as follows:\\
for $i\neq 0$,
\begin{equation*}
 v_{i,j}=\begin{blockarray}{*{2}{c}}
  \begin{block}{(c)c}
   {\bf 0}_2&\\
   \vdots&\\
   p_i& \leftarrow j\\
   \vdots &\\
   {\bf 0}_2&\\
   {\bf 0}_{3d+1}&\\
   {\bf 0}_{d-1}&\\
  \end{block}
 \end{blockarray}
 \in \Z^{2(d+1)+(3d+1)+(d-1)},
\end{equation*}
for $i=0$,
\begin{equation*}
  v_{0,j}=\begin{blockarray}{*{2}{c}}
  \begin{block}{(c)c}
   {\bf 0}_2&\\
   \vdots&\\
   p_0& \leftarrow j\\
   \vdots &\\
   {\bf 0}_2&\\
   -w_j&\\
   v'_j&\\
  \end{block}
 \end{blockarray}
 \in \Z^{2(d+1)+(3d+1)+(d-1)},
\end{equation*}
for $0\leq j\leq 3d$,
\begin{equation*}
  v_{3,j}=\begin{blockarray}{*{2}{c}}
  \begin{block}{(c)c}
   {\bf 0}_{2(d+1)}&\\
   e_j^{3d+1}&\\
   {\bf 0}_{d-1}&\\
  \end{block}
 \end{blockarray}
 \in \Z^{2(d+1)+(3d+1)+(d-1)},
\end{equation*}
and for $k=1,\dots,d-1$,
\begin{equation*}
 u_k=\begin{blockarray}{*{2}{c}}
  \begin{block}{(c)c}
   {\bf 0}_{2(d+1)}&\\
   {\bf 0}_{3d+1}&\\
   -e_k^{d-1}&\\
  \end{block}
 \end{blockarray}
 \in \Z^{2(d+1)+(3d+1)+(d-1)},\end{equation*}
where ${\bf 0}_\alpha$ is the zero vector in $\Z^\alpha$ and $e_k^\beta$ is the $k$-th standard basis of $\Z^\beta$.

\begin{defi}\label{definition:toric of P1113}
 Let
 \begin{align*}
  &P_0:=\{v_{0,0},v_{1,0},v_{2,0},v_{3,0},v_{3,1}\},\\
  &P_d:=\{v_{0,d},v_{1,d},v_{2,d},v_{3,3d-1},v_{3,3d}\},\\
  &P_i:=\{v_{0,i},v_{1,i},v_{2,i},v_{3,3i-1},v_{3,3i},v_{3,3i+1},u_i\}\; (1\leq i\leq d-1).
 \end{align*}
 Then, we define
 \begin{equation*}
  \Sigma_d
 \end{equation*}
 as a set of cones generated by the union of proper subsets (involving empty set) of $P_0,P_1,\dots,P_d$.
 (A cone corresponding to empty set is $\{0\}$).
\end{defi}
We have to show that $\Sigma_d$ is a fan.

\begin{theo}
 For arbitrary positive integer $d$, $\Sigma_d$ is a simplicial complete fan.
\end{theo}
In order to prove that $\Sigma_d$ is a simplicial complete fan, we should check the following claim:
\begin{lemm}\label{lemm:toric P1113}
 For all $v\in \R^{6d+2}$, there uniquely exist $a_{i,j}\in \R$ and $b_k\in\R$ that satisfy the following conditions.
 \begin{eqnarray*}
 \mbox{(a)} &&v=\sum_{i=0}^2\sum_{j=0}^d a_{i,j}v_{i,j} +\sum_{j=0}^{3d} a_{3,j}v_{3,j} + \sum_{k=1}^{d-1} b_ku_k,\\
 \mbox{(b)}&& {\rm min}(\{a_{0,0},a_{1,0},a_{2,0},a_{3,0},a_{3,1}\})=0,\\
   && {\rm min}(\{a_{0,d},a_{1,d},a_{2,d},a_{3,3d-1},a_{3,3d}\})=0,\\
 \mbox{(c)} &&{\rm min}(\{ a_{0,i},a_{1,i},a_{2,i},a_{3,3i-1},a_{3,3i},a_{3,3i+1},b_i \})=0\;\;(i=1,2,\cdots,d-1).
 \end{eqnarray*}
\end{lemm}
{\bf Proof.}
The following relations for $\{ v_{i,j}\},\{u_k\}$ hold:
\begin{align}
 &v_{0,0}+v_{1,0}+v_{2,0}+3v_{3,0}+2v_{3,1}+v_{3,2}-u_1=0,\\
 &v_{0,i}+v_{1,i}+v_{2,i}+v_{3,3i-2}+2v_{3,3i-1}+3v_{3,3i}+2v_{3,3i+1}+v_{3,3i+2}\nonumber \\
 &\qquad -u_{i-1}+2u_i-u_{i+1}=0, \quad (1\leq i \leq d-1)\\ 
 &v_{0,d}+v_{1,d}+v_{2,d}+v_{3,3d-2}+2v_{3,3d-1}+3v_{3,3d}-u_d=0.
\end{align}
We can easily show them by definition of $\Sigma_d$.

\vspace{0.5cm}
For all $v\in\R^{6d+2}$, it is clear that there uniquely exist real numbers $x_{i,j}$, $(i=1,2,\; 0\leq j\leq d)$, $x_{3,j}$, $(0\leq j\leq 3d)$, $y_k$, $(1\leq k\leq d-1)$ such that
\begin{equation*}
 v=\sum_{i=0}^d (x_{1,i}v_{1,i}+x_{2,i}v_{2,i})+\sum_{j=0}^{3d}x_{3,j}v_{3,j}+\sum_{k=1}^{d-1}y_ku_k.
\end{equation*}

Then, we obtain:
\begin{align*}
 &v\\
  =&v+\alpha_0(v_{0,0}+v_{1,0}+v_{2,0}+3v_{3,0}+2v_{3,1}+v_{3,2}-u_1)\\
  &+\alpha_d(v_{0,d}+v_{1,d}+v_{2,d}+v_{3,3d-2}+2v_{3,3d-1}+3v_{3,3d}-u_d)\\
  &+\sum_{i=1}^{d-1}\alpha_i(v_{0,i}+v_{1,i}+v_{2,i}+v_{3,3i-2}+2v_{3,3i-1}+3v_{3,3i}+2v_{3,3i+1}+v_{3,3i+2})\\
  =&(\alpha_0+x_{0,0})v_{0,0}+(\alpha_0+x_{1,0})v_{1,0}+(\alpha_0+x_{2,0})v_{2,0}+(3\alpha_0+x_{3,0})v_{3,0}\\
  &\qquad +(2\alpha_0+\alpha_1+x_{3,1})v_{3,1}\\
  &+(\alpha_d+x_{0,d})v_{0,d}+(\alpha_d+x_{1,d})v_{1,d}+(\alpha_d+x_{2,d})v_{2,d}+(3\alpha_d+x_{3,3d})v_{3,3d}\\
  &\qquad +(2\alpha_d+\alpha_{d-1}+x_{3,3d-1})v_{3,3d-1}\\
  &+\sum_{i=1}^{d-1}(\alpha_iv_{0,i}+(\alpha_i+x_{1,i})v_{1,i}+(\alpha_i+x_{2,i})v_{2,i}\\
  &\qquad+(\alpha_{i-1}+2\alpha_i+x_{3,3i-1})v_{3,3i-1}+(3\alpha_i+x_{3,3i})v_{3,3i}\\
  &\qquad+(2\alpha_i+\alpha_{i+1}+x_{3,3i+1})v_{3,3i+1}+(-\alpha_{i-1}+2\alpha_i-\alpha_{i+1}+y_i)u_i).
\end{align*}
Hence, we should show the following claim:
\begin{claim}
A map $F_d$ from $\R^{d+1}$ to $\R^{d+1}$:
\begin{align*}
 \left(
 \begin{array}{c}
 {\rm min}\{\alpha_0+z_0, 2\alpha_0+\alpha_1+x_{3,1}\}\\
 {\rm min}\{\alpha_1+z_1, \alpha_0+2\alpha_1+x_{3,2}, 2\alpha_1+\alpha_2+x_{3,4}, -\alpha_0+2\alpha_1-\alpha_2+y_1\}\\
 {\rm min}\{\alpha_2+z_2, \alpha_1+2\alpha_2+x_{3,5}, 2\alpha_2+\alpha_3+x_{3,7}, -\alpha_1+2\alpha_2-\alpha_3+y_2\}\\
 \vdots\\
 {\rm min}\{\alpha_{d-1}+z_{d-1}, \alpha_{d-2}+2\alpha_{d-1}+x_{3,3d-4}, 2\alpha_{d-1}+\alpha_d+x_{3,3d-2}, \\
 \qquad -\alpha_{d-2}+2\alpha_{d-1}-\alpha_d+y_{d-1}\}\\
 {\rm min}\{\alpha_d+z_d, \alpha_{d-1}+2\alpha_d+x_{3,3d-1}\}
 \end{array}
 \right)
\end{align*}
is bijective, where $z_i:=\max\{0,x_{0,i},x_{1,i},x_{2,i},x_{3,3i}/3\}$.
\end{claim}

The following two lemmas lead us to proof of the above claim:
\begin{lemm}\label{lemm:coherently oriented}
$F_d(\alpha)$ is coherently oriented piecewise affine map (i.e. for any component $P\subset\R^{d+1}$ which $F(\alpha)$ is linear, ${\rm det}(F|_P)$ is positive).
\end{lemm}

\begin{lemm}\label{lemm:recession function}
 The {\rm recession function} of $F_d$
 \begin{align*}
  &F_d^\infty(\alpha):=\\
  &\left(
  \begin{array}{c}
  {\rm min}\{\alpha_0, 2\alpha_0+\alpha_1\}\\
  {\rm min}\{\alpha_1, \alpha_0+2\alpha_1, 2\alpha_1+\alpha_2, -\alpha_0+2\alpha_1-\alpha_2\}\\
  {\rm min}\{\alpha_2, \alpha_1+2\alpha_2, 2\alpha_2+\alpha_3, -\alpha_1+2\alpha_2-\alpha_3\}\\
  \vdots\\
  {\rm min}\{\alpha_{d-1}, \alpha_{d-2}+2\alpha_{d-1}, 2\alpha_{d-1}+\alpha_d, -\alpha_{d-2}+2\alpha_{d-1}-\alpha_d\}\\
  {\rm min}\{\alpha_d, \alpha_{d-1}+2\alpha_d\}
  \end{array}
  \right)
 \end{align*}
 is bijective.
\end{lemm}
If these lemmas are proved, we can use the following theorem:
\begin{theo}[Theorem 2.5.1. of \cite{Sch}]
 A coherently oriented piecewise affine function is a homeomorphism if and only if its recession function is a homeomorphism.
\end{theo}

Hence, $F_d(\alpha)$ is also a homeomorphism (i.e., bijective).

\subsection{Proof of Lemma \ref{lemm:coherently oriented}.}
 The coefficient matrix of $F_d$ is given as follows.

\vspace{1cm}
\scalebox{0.7}{
\begin{minipage}{\linewidth}
\begin{align*}\label{coefficient matrix}
 \left(
 \begin{array}{cccccc:c:ccccc:c:c:ccccc}
 1+\epsilon_0&\epsilon_0&0&\cdots&0&0&0&&&&&&&&&&&\\
 \epsilon_1&2&\epsilon_1'&\cdots&0&0&0&&&&&&&&&&&\\
  0&\epsilon_2&2&\cdots&0&0&0&&&&&&&&&&&\\
  \vdots&\vdots&\vdots&\ddots&\vdots&\vdots&\vdots&&&&&&&&&&&\\
  0&0&0&\cdots&2&\epsilon_{k_1-2}'&0&&&&&&&&&&&\\ 
  0&0&0&\cdots&\epsilon_{k_1-1}&2&\epsilon_{k_1-1}'&&&&&&&&&&&\\ \hdashline
 &&&&&&1&&&&&&&&&&&\\ \hdashline
 &&&&&&\epsilon_{k_1+1}&2&\epsilon_{k_1+1}'&\cdots&0&0&&&&&&&\\
 &&&&&&0&\epsilon_{k_1+2}&2&\cdots&0&0&&&&&&&\\
 &&&&&&\vdots&\vdots&\vdots&\ddots&\vdots&\vdots&&&&&&&\\
 &&&&&&0&0&0&\cdots&2&\epsilon_{k_2-2}'&&&&&&&\\
 &&&&&&0&0&0&\cdots&\epsilon_{k_2-1}&2&\epsilon_{k_2-1}'&&&&&\\ \hdashline
 &&&&&&&&&&&&1&&&&&\\ \hdashline
 &&&&&&&&&&&&&\ddots&&&&&\\ \hdashline
 &&&&&&&&&&&&&&2&\epsilon_{k_r+1}'&\cdots&0&0\\
 &&&&&&&&&&&&&&\epsilon_{k_r+2}&2&\cdots&0&0\\
 &&&&&&&&&&&&&&\vdots&\vdots&\ddots&\vdots&\vdots\\
 &&&&&&&&&&&&&&0&0&\cdots&2&\epsilon_{d-1}'\\
 &&&&&&&&&&&&&&0&0&\cdots&\epsilon_d&1+\epsilon_d
 \end{array}
 \right)
\end{align*}
\end{minipage}}
\vspace{1cm}

where $0<k_1<k_2<\cdots<k_r=d$ and
$(\epsilon_0,\epsilon_1,\epsilon_1',\dots,\epsilon_{d-1},\epsilon_{d-1}',\epsilon_d)\in \{0,\pm1\}^{2d}$.

In the above matrix, we can assume that for $1\leq i\leq d-1$,
\begin{equation*}
 (\epsilon_i,\epsilon_i')=
 \begin{cases}
  (1,0),\; (0,1),\; {\rm or}\; (-1,-1),&(i\notin K)\\
  (0,0), &(i\in K)
 \end{cases}
\end{equation*}
and
\begin{equation*}
 (\epsilon_0,\epsilon_d)\in\{0,1\}^2,
\end{equation*}
where $K:=\{k_1,k_2,\dots,k_r\}$.

Then what we have to show is that determinant of the matrix is positive.
From a formula
\begin{equation*}
 {\rm det}\left(
 \begin{array}{cc}
  X&0\\
  Z&Y
 \end{array}
 \right)={\rm det}(X){\rm det}(Y),
\end{equation*}
and elementary operations of matrix, we can reduce the problem to check positivity 
of determinant of the following matrix:
\begin{equation*}
 B_k:=\left(
 \begin{array}{cccccc}
 2&1&0&\cdots&0&0\\
 -1&2&-1&\cdots&0&0\\
  0&-1&2&\cdots&0&0\\
  \vdots&\vdots&\vdots&\ddots&\vdots&\vdots\\
  0&0&0&\cdots&2&-1\\ 
  0&0&0&\cdots&1&2
 \end{array}
 \right)\in M_{k+1,k+1},
\end{equation*}
where $k>0$.
We can easily compute its determinant:
\begin{equation*}
 {\rm det}B_k=9k-6>0.
\end{equation*}
$\Box$

\subsection{Proof of Lemma \ref{lemm:recession function}.}
It is clear that
\begin{align*}
 F_d^\infty(t\alpha)=tF_d^\infty(\alpha),
\end{align*}
for all $t\in\R_{\geq 0}$
Hence, in order to prove that $F_d^\infty$ is injective, we only have to show that
\begin{align*}
 G:=\pi\circ F_d^\infty|_{S^d}:S^d\rightarrow S^d,
\end{align*}
is injective.
where $\pi:\R^{d+1}-\{0\}\rightarrow S^d$ is a projection.

Obviously $G$ is continuous.
Let $\tilde{G}$ be a smoothing of $G$. Note that we can take volume of smoothing locus of $G$ as small as we like
because non-smooth locus of $G$ is measure $0$.  
Then, $G(\alpha)\neq -\alpha$ holds for all $\alpha\in S^d$ since the diagonal elements of $F_d^\infty$ are all positive.
Therefore,
\begin{equation*}
 H(t,\alpha):=\pi(t\alpha+(1-t)\tilde{G}(\alpha))
\end{equation*}
gives us a homotopy from $\tilde{G}$ to the identity mapping $id_{S^d}$.
Thus, the mapping degree of $\tilde{G}$ is $1$.
Hence, $\tilde{G}$ is injective since Jacobian of $\tilde{G}$ is always positive.
Accordingly, $G$ is also injective.  $\Box$

\subsection{Some Properties of $\widetilde{Mp}_{0,2}(\bP(1,1,1,3),d)$.}
Since the fan $\Sigma_d$ which we constructed in the previous section is complete and simplicial, the corresponding toric variety $X_{\Sigma_d}$ is a compact orbifold.

In this subsection, we will show that the toric variety $X_{\Sigma_d}$ realizes the action (\ref{eq:action P1113}), and prove that $\widetilde{Mp}_{0,2}(\bP(1,1,1,3),d)$ is a compact toric orbifold.

First, we determine the primitive collections of the fan $\Sigma_d$.
\begin{lemm}
 The primitive collections of the fan $\Sigma_d$ are
 \begin{align*}
  &P_0:=\{v_{0,0},v_{1,0},v_{2,0},v_{3,0},v_{3,1}\},\\
  &P_d:=\{v_{0,d},v_{1,d},v_{2,d},v_{3,3d-1},v_{3,3d}\},\\
  &P_i:=\{v_{0,i},v_{1,i},v_{2,i},v_{3,3i-1},v_{3,3i},v_{3,3i+1},u_i\}\; (1\leq i\leq d-1).
 \end{align*}
\end{lemm}
{\bf Proof.}
By Definition \ref{definition:toric of P1113}, it is clear that they are primitive collections of $\Sigma_d$.
If $P$ is a primitive collection of $\Sigma_d$, then $P$ does not generate any cone of $\Sigma_d$.
Hence $P$ has to contain some $P_i$ $(i=0,1,\dots,d)$.
If $P_i$ is proper subset of $P$ for some $i=0,1,\dots,d$, then $P_i$ does not generate any cone of $\Sigma_d$.
Therefore, $P$ should not be a primitive collection.
Accordingly, $P=P_i$. $\Box$

We introduce the following notation,
\begin{align*}
 \C^{\Sigma_d(1)}=\{ ({\bf a}_0,{\bf a}_1,\dots,{\bf a}_d,a_{3,0},a_{3,1},\dots,a_{3,3d},u_1,u_2,\dots,u_{d-1}) \; | \; {\bf a}_i \in \C^3,\; a_{3,j}, u_i\in \C \}.
\end{align*}
where $\Sigma_d(1)$ is a collection of 1-dimensional cones of $\Sigma_d$.
Then, we define a subset $Z(\Sigma_d)$ of $\C^{\Sigma_d(1)}$ as follows.
\begin{equation*}
 Z(\Sigma_d)=\left\{ x\in \C^{\Sigma_d(1)}\;|\; 
\begin{array}{ll}
 ({\bf a}_0,a_{3,0},a_{3,1})=0,&\\
 ({\bf a}_i,a_{3,3i-1},a_{3,3i},a_{3,3i+1},u_i)=0, & (1\leq i\leq d-1)\\
 ({\bf a}_d,a_{3,3d-1},a_{3,3d})=0,
\end{array}
\right\}.
\end{equation*}
Note that the toric variety corresponding to the fan $\Sigma_d$ is given by quotient space $(\C^{\Sigma_d(1)}\backslash Z(\Sigma_d))/G$, where $G:={\rm Hom}_\Z(A_{{\rm dim}(X_{\Sigma_d})-1}(X_{\Sigma_d}),\C^*)$ (see \cite{Cox} or Chap. 3 of \cite{CK}).
The $G$-action is determined as follows.

Let $[D_{i,j}]$ (resp. $[U_k]$) be a divisor class that corresponds to 1-dimensional cone $v_{i,j}$ (resp. $u_k$).
Furthermore, let
\begin{equation*}
 n:={\rm dim}(X_{\Sigma_d})=6d+2.
\end{equation*}
Recall the following exact sequence:
\begin{equation}\label{eq:ex seq of toric}
 0\rightarrow M\rightarrow \Z^{\Sigma_d(1)}\rightarrow A_{n-1}(X_{\Sigma_d})\rightarrow 0.
\end{equation}
Here $M=\Z^{6d+2}$.
Note that $A_{n-1}(X_{\Sigma_d})$ is generated by $[D_{i,j}]$ and $[U_k]$.
$M\rightarrow \Z^{\Sigma_d(1)}$ and $\Z^{\Sigma_d(1)}\rightarrow A_{n-1}(X_{\Sigma_d})$ are given by
\begin{align*}
 &M\rightarrow \Z^{\Sigma_d(1)};m\mapsto (\left<m,v_\rho\right>)_{\rho\in\Sigma_d(1)}\\
 &\Z^{\Sigma_d(1)}\rightarrow A_{n-1}(X_{\Sigma_d});(a_\rho)_{\rho\in\Sigma_d(1)}\mapsto \sum_{\rho\in\Sigma_d(1)}a_\rho[D_\rho].
\end{align*}
Here $v_\rho\in M$ is a generator of 1-dimensional cone $\rho\in\Sigma_d(1)$ and $D_\rho$ is a divisor corresponding to $\rho\in\Sigma_d(1)$.

By the exact sequence (\ref{eq:ex seq of toric}) and definition of $v_{i,j}$ and $u_k$, we obtain the following relations of Chow group $A_{n-1}(X_{\Sigma_d})$:
\begin{equation}\label{eq:rel for div}
 \begin{cases}
  [D_{0,j}]=[D_{1,j}]=[D_{2,j}]\, (0\leq j\leq d),\\
  [D_{3,3j}]=3[D_{0,j}],\, (0\leq j\leq d)\\
  [D_{3,3j+1}]=2[D_{0,j}]+[D_{0,j+1}],\, (0\leq j\leq d-1)\\
  [D_{3,3j+2}]=[D_{0,j}]+2[D_{0,j+1}],\, (0\leq j\leq d-1)\\
  [U_k]=-[D_{0,k-1}]+2[D_{0,k}]-[D_{0,k+1}]\, (1\leq k\leq d-1).
 \end{cases}
\end{equation}
From these relations, it is easily shown that $G={\rm Hom}_{\Z}(A_{n-1}(X_{\Sigma_d}),(\C^*))\cong (\C^*)^{d+1}$.
Let $\lambda_i:=g([D_{i,1}])$ ($g\in G={\rm Hom}_{\Z}(A_{n-1}(X_{\Sigma_d}),(\C^*))$).
The above relation tells us that $g([U_k])=\lambda_{k-1}^{-1}\lambda_k^2\lambda_{k+1}^{-1}$, and the $({\C^*})^{d+1}$-action turns out to be,
\begin{align}
 &(\lambda_0,\dots,\lambda_d)\cdot ({\bf a}_0,{\bf a}_1,\dots,{\bf a}_d,a_{3,0},a_{3,1},\dots,a_{3,3d},u_1,u_2,\dots,u_{d-1})\nonumber \\
 =&(\lambda_0{\bf a}_0,\lambda_1{\bf a}_1,\dots,\lambda_d{\bf a}_d,\lambda_0^3a_{3,0},\lambda_0^2\lambda_1a_{3,1},\lambda_0\lambda_1^2a_{3,2},\lambda_1^3a_{3,3},\lambda_1^2\lambda_2a_{3,4},\dots,\lambda_d^3a_{3,3d},\nonumber\\
&\lambda_0^{-1}\lambda_1^2\lambda_2^{-1}u_1,\lambda_1^{-1}\lambda_2^2\lambda_3^{-1}u_2,\dots,\lambda_{d-2}^{-1}\lambda_{d-1}^2\lambda_d^{-1}u_{d-1}). \nonumber
\end{align}

When $u_k=1$ for all $k=1,2,\dots,d-1$, by setting $\lambda_i=(\mu')^{d-i}(\nu')^i$, we obtain the action which is similar to (\ref{eq:action P1113}).

Accordingly, we can identify $X_{\Sigma_d}=\widetilde{Mp}_{0,2}(\bP(1,1,1,3),d)$.
Hence, we obtain Theorem \ref{theo:properties MpP1113}.

\subsection{The Chow Ring of $\widetilde{Mp}_{0,2}(\bP(1,1,1,3),d)$ (proof of proposition \ref{prop:Chow ring of MpP1113}).}
In this subsection, we will compute the Chow ring of $\widetilde{Mp}_{0,2}(\bP(1,1,1,3),d)$.
The recipe of computation is the same as the one in \cite{S}.

Let us recall the structure of Chow ring of general toric variety is given by
\begin{equation*}
 A^*(X_\Sigma)\cong\C[x_\rho|\rho\in\Sigma(1)]/(P(\Sigma)+SR(\Sigma)).
\end{equation*}
Here, 
\begin{align*}
 P(\Sigma)&:=\left<\sum_{\rho\in\Sigma(1)}\left<m,v_\rho\right>x_\rho|m\in M\right>\\
 SR(\Sigma)&:=\left<x_{\rho_1}\cdots x_{\rho_k}|\{\rho_1,\dots,\rho_k\} {\rm \,is\, a\, primitive\, collection\, of\,} \Sigma\right>.
\end{align*}
(see \cite{Ful}).

{\bf Proof of proposition \ref{prop:Chow ring of MpP1113}.}
It is easily see that
\begin{equation*}
 \C[x_\rho|\rho\in\Sigma_d(1)]/P(\Sigma_d)\cong \C[H_0,H_1,\dots,H_d]
\end{equation*}
by setting $H_j:=[D_{0,j}]$ $(j=0,1,\dots,d)$ and the relations (\ref{eq:rel for div}).

Recall that the primitive sets of $\Sigma_d$ are
\begin{align*}
 &P_0:=\{v_{0,0},v_{1,0},v_{2,0},v_{3,0},v_{3,1}\},\\
 &P_d:=\{v_{0,d},v_{1,d},v_{2,d},v_{3,3d-1},v_{3,3d}\},\\
 &P_i:=\{v_{0,i},v_{1,i},v_{2,i},v_{3,3i-1},v_{3,3i},v_{3,3i+1},u_i\}\; (1\leq i\leq d-1).
\end{align*}
Then, the Stanley-Reisner ideal is
\begin{align*}
 SR(\Sigma_d)=&(H_0^4(2H_0+H_1),H_1^4(H_0+2H_1)(2H_1+H_2)(-H_0+2H_1-H_2),\\
 &\dots,H_{d-1}^4(H_{d-2}+2H_{d-1})(2H_{d-1}+H_d)(-H_{d-2}+2H_{d-1}-H_d),\\
 &H_d^4(H_{d-1}+2H_d)).
\end{align*}
$\Box$

\subsection{The Intersection Numbers $w({\cal O}_{h^a}{\cal O}_{h^b})_{0,d}$.}
We use
\begin{align*}
 {\rm Vol}_d:&=\left(\prod_{i=0}^d[D_{1,i}][D_{2,i}]\right)\cdot\left(\prod_{j=0}^{3d}[D_{3,j}]\right)\cdot\left(\prod_{k=1}^{d-1}[U_k]\right)\\
 &=3^{d+1}\left(\prod_{i=0}^d H_i^3\right)\cdot\left(\prod_{i=0}^{d-1}(2H_i+H_{i+1})(H_i+2H_{i+1})\right)\\
 &\quad \times\left(\prod_{k=1}^{d-1}(-H_{k-1}+2H_k-H_{k+1})\right).
\end{align*}
as a volume form of $A^*(\widetilde{Mp}_{0,2}(\bP(1,1,1,3),d))$ since it corresponds to Poincar\'{e} dual of a smooth point on the toric variety.

Let us explain Definition \ref{definition:intersection num of MpP1113} of intersection numbers of $\widetilde{Mp}_{0,2}(\bP(1,1,1,3),d)$ again:
\begin{defi}
 Let
 \begin{equation*}
  e^{6}(x,y):=\prod_{j=0}^6 ((6-j)x+y).
 \end{equation*}
 Then, we define the intersection number $w({\cal O}_a{\cal O}_b)_{0,d}$ over $\widetilde{Mp}_{0,2}(\bP(1,1,1,3),d)$ as
 \begin{eqnarray*}
  w({\cal O}_{z^a}{\cal O}_{z^b})_{0,d}:=\int_{\widetilde{Mp}_{0,2}(\bP(1,1,1,3),d)} H_0^aH_d^b \cdot \frac{\displaystyle{\prod_{i=1}^de^{6}(H_{i-1},H_i)}}{\displaystyle{\prod_{i=1}^{d-1}6H_i}}.
 \end{eqnarray*}
\end{defi}
This intersection number is a quasi map analogue of genus $0$ degree $d$ two point GW invariants of K3 surface in $\bP(1,1,1,3)$.
\begin{prop}
For $\Omega\in A^*(\widetilde{Mp}_{0,2}(\bP(1,1,1,3),d))$, the following equality holds. 
\begin{equation*}
 \int_{\widetilde{Mp}_{0,2}(\bP(1,1,1,3),d)}\Omega=\prod_{i=0}^d\left(\frac1{2\pi\sqrt{-1}}\oint_{C_j}{dz_j}\right) \frac{\tilde{\Omega}}{R}
\end{equation*}
where
\begin{eqnarray*} 
R&=&3^{d+1}\biggl(\prod_{i=0}^d(z_i)^4\biggr)(2z_0+z_1)\biggl(\prod_{i=1}^{d-1}(z_{i-1}+2z_i)(2z_i+z_{i+1})\biggr)(z_{d-1}+2z_d)\times\nonumber\\
&&\biggl(\prod_{i=1}^{d-1}(2z_j-z_{j-1}-z_{j+1})\biggr). 
\end{eqnarray*}
$\frac{1}{2\pi\sqrt{-1}}\oint_{C_j}\;(j=1,2,\cdots,d-1)$ means taking residues at $z_j=0$, $z_{j}=-\frac{z_{j-1}}{2}$ $z_{j}=-\frac{z_{j+1}}{2}$, $\frac{z_{j-1}+z_{j+1}}{2}$ and $\frac{1}{2\pi\sqrt{-1}}\oint_{C_0}$ (resp. $\frac{1}{2\pi\sqrt{-1}}\oint_{C_d}$) means
taking residues at $z_{0}=0$, $z_{0}=-\frac{z_{1}}{2}$ (resp. $z_{d}=0$, $z_{d}=-\frac{z_{d-1}}{2}$).
$\tilde{\Omega}$ is a polynomial obtained from turning $H_i$ into $z_i$ in $\Omega$.
\label{residue}
\end{prop}
{\bf Proof.} Obviouly, the correspondence $Res:\tilde{\Omega}\rightarrow \prod_{i=0}^d\left(\frac1{2\pi\sqrt{-1}}\oint_{C_j}{dz_j}\right)
 \frac{\tilde{\Omega}}{R}$ defines a linear map from $\mathbb{C}[z_{0},z_{1},\cdots,z_{d}]$ to $\mathbb{C}$. Let ${\cal I}$ be 
an ideal of $\mathbb{C}[z_{0},z_{1},\cdots,z_{d}]$ generated by, 
\begin{eqnarray}
&&r_{0}:=z_{0}^4(2z_{0}+z_{1}), \nonumber\\
&&r_{1}:=z_{1}^{4}(2z_{1}+z_{0})(2z_{1}+z_{2})(2z_{1}-z_{0}-z_{2}),\nonumber\\
&&\hspace{1cm}\vdots\nonumber\\
&&r_{d-1}:=z_{d-1}^{4}(2z_{d-1}+z_{d-2})(2z_{d-1}+z_{d})(2z_{d-1}-z_{d-2}-z_{d}),\nonumber\\
&&r_{d}:=z_{d}^4(2z_{d}+z_{d-1}).
\end{eqnarray}
If $\tilde{\Omega}$ takes the form $r_{i}\cdot f\;(f\in \mathbb{C}[z_{0},z_{1},\cdots,z_{d}])$, $Res(r_{i}\cdot f)=0$ since 
the integrand is holomorphic at the points where we take residues of $z_{i}$. 
Moreover, we can easily see by degree counting that $Res(\tilde{\Omega})$ vanishes if $\tilde{\Omega}$ is a homogeneous polynomial whose degree is not equal to $6d+2$.    
Therefore, $Res$ turns out to be a map from $\bigl(\mathbb{C}[z_{0},z_{1},\cdots,z_{d}]/{\cal I}\bigr)_{6d+2}$ (homogeneous degree 
$6d+2$ part of  $\mathbb{C}[z_{0},z_{1},\cdots,z_{d}]/{\cal I}$) to $\mathbb{C}$. 
By investigating Hilbert polynomial of the ring $\mathbb{C}[z_{0},z_{1},\cdots,z_{d}]/{\cal I}$ (which is isomorphic to 
$A^*(\widetilde{Mp}_{0,2}(\bP(1,1,1,3),d))$), $\bigl(\mathbb{C}[z_{0},z_{1},\cdots,z_{d}]/{\cal I}\bigr)_{6d+2}$ turns 
out to be $1$-dimensional. Hence we only have to check the following equality:
\begin{eqnarray}
&&Res(\widetilde{\rm Vol}_{d})=1,\nonumber\\
&\Longleftrightarrow& \prod_{i=0}^d\left(\frac1{2\pi\sqrt{-1}}\oint_{C_j}{dz_j}\right) \biggl(\prod_{j=0}^{d}\frac{1}{z_{j}}\biggr)=1. 
\end{eqnarray}  
But this is obvious.            $\Box$.              
                       
\begin{rema}
We can easily see from the proof of Proposition \ref{residue} that result of the residue integral does not depend on 
order of integration with respect to subscript $j$ of $z_{j}$.  
\end{rema}

\section{Proof of Main Theorem \ref{theorem:main}}\label{section:j-inv}

In this section, we prove our main theorem of this paper.
Let us restate the theorem here.
\begin{theo}
 Let $w_d$ be the $d$-th expansion coefficient of inverse function of $-\log(j(\tau))$.
 Then
 \begin{equation*}
  w_d=\frac12w({\cal O}_{z^1}{\cal O}_{z^0})_{0,d}.
 \end{equation*}
\end{theo}
In order to prove it, we should show that
\begin{equation}\label{eq:goal}
 \frac{f_1(e^x)}{f_0(e^x)} = x + \sum_{d=1}^\infty w_d e^{dx} = x + \sum_{d=1}^\infty \frac12 w({\cal O}_{z^1}{\cal O}_{z^0})_{0,d}e^{dx},
\end{equation}
where
\begin{align*}
 f_0(z)&:=\sum_{d=0}^\infty \frac{2^{3d}\cdot (6d-1)!!}{(d!)^3}z^d,\\
 f_1(z)&:=f_0(z)(\log(z))+\sum_{d=0}^\infty (\sum_{j=1}^{3d}\frac{6}{2j-1}-\sum_{k=1}^{d}\frac{3}{j})\frac{(6d-1)!!}{(d!)^3}z^d
\end{align*}
are two solutions of Picard-Fuchs equation (\ref{eq:P-F eq of K3}) as subsection \ref{subsec:PF eq for j-inv}.

Let us introduce the following generating functions.
\begin{align*}
 L_0(e^x)&:=1+\sum_{d=1}^\infty \frac{d}{2}w({\cal O}_{z^2}{\cal O}_{z^{-1}})_{0,d}e^{dx},\\
 L_1(e^x)&:=1+\sum_{d=1}^\infty \frac{d}{2}w({\cal O}_{z^1}{\cal O}_{z^0})_{0,d}e^{dx}.
\end{align*}

First lemma claims that $L_0(e^x)$ gives us the solution $f_0(e^x)$ of Picard-Fuchs equation (\ref{eq:P-F eq of K3}).
\begin{lemm}\label{lemm:hol sol}
 \begin{equation*}
  f_0(e^x)=L_0(e^x).
 \end{equation*}
\end{lemm}
{\bf Proof.}
We have to prove that
\begin{eqnarray}
&& \frac{2^{3d}(6d-1)!!}{(d!)^3}=\frac{d}{2}w({\cal O}_{z^2}{\cal O}_{z^{-1}})_{0,d}\nonumber\\
&\Longleftrightarrow& \frac{1}{2}w({\cal O}_{z^2}{\cal O}_{z^{-1}})_{0,d}= \frac{1}{d}\cdot\frac{2^{3d}(6d-1)!!}{(d!)^3} 
\label{as}
\end{eqnarray}
for all $d$.
From Proposition \ref{residue},
we have the following equality:
\begin{eqnarray}
&&\frac{1}{2}w({\cal O}_{z^2}{\cal O}_{z^{-1}})_{0,d}\nonumber\\
&=&\frac{1}{2}\left(\frac{1}{2\pi\sqrt{-1}}\oint_{C_j}dz_j\right)z_{0}^{2}\left(\prod_{j=0}^{d-1}e^{6}(z_{j},z_{j+1})\right)
\left(\prod_{j=1}^{d-1}\frac{1}{6z_{j}}\right)\frac{1}{z_{d}}\cdot\frac{1}{R}\nonumber\\
&=&\frac{1}{2\cdot 3^{d+1}}\prod_{j=0}^d\left(\frac1{2\pi\sqrt{-1}}\oint_{C_j}\frac{dz_j}{z_j^4}\right)z_0^2\left(\prod_{j=1}^d\tilde{e}(z_{j-1},z_j)\right)\nonumber\\
&&\quad\times \left(\prod_{j=1}^{d-1}\frac1{6z_j(2z_j-z_{j-1}-z_{j+1})}\right)\frac1{z_d},
\label{00}
\end{eqnarray}
where
\begin{equation*}
 \tilde{e}(x,y):=2^4\cdot 3^2xy\prod_{i=0}^2((2i+1)x+(5-2i)y).
\end{equation*}
Since the integral is holomorphic at $z_{0}=-\frac{z_{1}}{2},\;z_{j}=-\frac{z_{j-1}}{2}, -\frac{z_{j+1}}{2}\;(j=1,\cdots,d-1),\; z_{d}=-\frac{z_{d-1}}{2}$, 
we only have to take residues at $z_{0}=0, z_{j}=0,\frac{z_{j-1}+z_{j+1}}{2}\;(j=1,\cdots,d-1), z_{d}=0$.
Note that the integral does not depend on order of integration.
Therefore, we integrate the last line of (\ref{00}) in ascending order of subscript $j$.

First, we integrate out $z_0$ variable.
By picking up the factors containing $z_0$, integration is done as follows:
\begin{align*}
 &\frac1{2\cdot 3^{d+1}}\frac1{2\pi\sqrt{-1}}\oint_{C_{(0)}}\frac{dz_0}{z_0}2^4\cdot 3^2z_1\prod_{i=0}^2((2i+1)z_0+(5-2i)z_1)\frac1{2z_1-z_0-z_2}\\
 =&\frac{2^3\cdot 5!!}{3^{d-1}}\frac{z_1^4}{2z_1-z_2}.
\end{align*}
Then we integrate $z_1$ variable.
Since the integrand is holomorphic at $z_1=0$, we only have to take residue at $z_1=z_2/2$.
\begin{align*}
 &\frac1{2\cdot 3^{d+1}}\prod_{j=0}^d\frac1{2\pi\sqrt{-1}}\oint_{C_{z_2/2}}dz_1\cdot 2^3\cdot 3 z_2\prod_{i=0}^2((2i+1)z_1+(5-2i)z_2)\frac1{2z_2-z_1-z_3}\\
 =&\frac12\cdot \frac{2^6\cdot 11!!}{3^{d-2}\cdot (2!)^3}\cdot\frac{z_2^4}{\frac32z_2-z_2}.
\end{align*}
Here, we use the identity:
\begin{equation*}
 \prod_{i=0}^2((2i+1)\frac{z_2}{2}+(5-2i)z_2)=(z_2)^3\prod_{i=0}^2\frac{11-2i}{2}=z_2^3\frac{11!!}{5!!\cdot 2^3}.
\end{equation*}
Integration of $z_i$ ($i=1,2,\dots,d-1$) goes in the same way.
We only have to take residue at $z_j=\frac{j}{j+1}z_{j+1}$.
After finishing integration of $z_{d-1}$, what remains to do is the following integration.
\begin{align*}
 \frac1{d}\cdot\frac{2^{3d}\cdot(6d-1)!!}{(d!)^3}\frac1{2\pi\sqrt{-1}}\oint_{C_{(0)}}\frac{dz_d}{z_d}.
\end{align*}
Hence we obtain the equality (\ref{as}). $\Box$

The following is the second lemma:
\begin{lemm}\label{lemm:second}
 \begin{align*}
  &\frac1{2\cdot 3^{d+1}}\prod_{j=0}^d\left(\frac1{2\pi\sqrt{-1}}\oint_{C_j}\frac{dz_j}{z_j^4}\right)z_0z_1\left(\prod_{j=1}^d\tilde{e}(z_{j-1},z_j)\right)\\
  &\quad\times \left(\prod_{j=1}^{d-1}\frac1{6z_j(2z_j-z_{j-1}-z_{j+1})}\right)\frac1{z_d}\\
  .&=\frac1d\cdot \frac{2^{3d}(6d-1)!!}{(d!)^3}\left(1-\frac1d+\sum_{j=1}^{3d}\frac{6}{2j-1}-\sum_{j=1}^d\frac3j\right).
 \end{align*}
\end{lemm}
{\bf Proof.}
In the same way as the proof previous lemma, we begin by integrating out $z_0$ variable.
\begin{align}\label{eq:0.8}
 &\frac1{2\cdot 3^{d+1}}\frac1{2\pi\sqrt{-1}}\oint_{C_{(0)}}\frac{dz_0}{z_0^2}2^4\cdot 3^2z_1\prod_{i=0}^2((2i+1)z_0+(5-2i)z_1)\frac1{2z_1-z_0-z_2}\nonumber\\
 =&\frac{2^3\cdot 5!!}{3^{d-1}}\left((a_1)\frac{z_1^3}{2z_1-z_2}+\frac{z_1^4}{(2z_1-z_2)^2}\right),
\end{align}
where
\begin{align*}
 a_1:=\sum_{i=0}^2\frac{2i+1}{5-2i}.
\end{align*}
In deriving (\ref{eq:0.8}), we used the following equality:
\begin{align*}
 &\frac{\partial}{\partial z_0}\left(\prod_{i=0}^2((2i+1)z_0+(5-2i)z_1)\frac1{2z_1-z_0-z_2}\right)\\
 =&\left(\prod_{i=0}^2((2i+1)z_0+(5-2i)z_1)\right)\\
 &\quad \times \frac1{2z_1-z_0-z_2}\left(\sum_{i=0}^2\frac{2i+1}{(2i+1)z_0+(5-2i)z_1}+\frac1{2z_1-z_0-z_2}\right).
\end{align*}
Since we have another $z_1$ factor in the integrand, it becomes holomorphic at $z_1=0$ after integration of $z_0$.
Hence integration of $z_1$ variable is done by taking residue at $z_1=z_2/2$.
\begin{align}\label{eq:0.11}
 &\frac{3^3\cdot 5!!}{3^{d-1}}\frac1{2\pi\sqrt{-1}}\oint_{C_{(z_2/2)}}dz_1(a_1\frac1{2z_1-z_2}+\frac{z_1}{(2z_1-z_2)^2})\nonumber \\
 &\quad \times 2^4\cdot 3^2z_2\prod_{i=0}^2((2i+1)z_1+(5-2i)z_2)\frac1{2z_1-z_0-z_2}\nonumber \\
 =&\frac12\cdot \frac{2^6\cdot 11!!}{3^{d-2}(2!)^3}\left(a_2\frac{z_2^4}{\frac32z_2-z_3}+\frac14\frac{z_2^5}{(\frac32 z_2-z_3)^2}\right).
\end{align}
Here, $a_2$ is given by
\begin{equation*}
 a_2:=a_1+\frac{1}{2\cdot 1}+\frac14\sum_{i=0}^2\frac{4i+2}{11-2i}.
\end{equation*}
In deriving (\ref{eq:0.11}), we used the following equality:
\begin{align*}
 &\frac{\partial}{\partial z_1}\left(z_1(\prod_{i=0}^2((2i+1)z_1+(5-2i)z_2))\frac1{2z_2-z_1-z_3}\right)\\
 =&z_1\left(\prod_{i=0}^2((2i+1)z_1+(5-2i)z_2)\right)\frac1{2z_2-z_1-z_3}\\
  &\quad \times\left(\frac1{z_1}+\sum_{i=0}^2\frac{2i+1}{(2i+1)z_1+(5-2i)z_2}+\frac1{2z_2-z_1-z_3}\right).
\end{align*}
Integration of $z_j$ ($j=2,3,\dots,d-2$) goes in the same way. After finishing integration of $z_{d-1}$, the LHS of this lemma becomes
\begin{equation*}
 \frac1d\cdot\frac{2^3d\cdot(6d-1)!!}{(d!)^3}a_d\frac1{2\pi\sqrt{-1}}\oint_{C_{(0)}}\frac{dz_d}{z_d},
\end{equation*}
where
\begin{align*}
 a_d&=\sum_{j=2}^d\frac1{j(j-1)}+\sum_{j=1}^d\frac{1}{j^2}\sum_{i=0}^2\frac{j(2i+1)}{6j-1-2i}\\
 &=\sum_{j=2}^d\left(\frac1{j-1}-\frac1j\right)+\sum_{j=1}^d\sum_{i=0}^2\frac{(2i+1)}{(6j-1-2i)j}\\
 &=1-\frac1d+\sum_{j=1}^d\sum_{i=0}^2(\frac6{6j-1-2i}-\frac6{6j})\\
 &=1-\frac1d+\sum_{j=1}^{3d}\frac6{2j-1}-\sum_{j=1}^d\frac3j.
\end{align*}
Integration of $z_d$ immediately leads us to the assertion of the lemma. $\Box$

We prove the last lemma:
\begin{lemm}\label{lemm:last}
 \begin{align}
  &\frac12 w({\cal O}_{z^1}{\cal O}_{z^0})_{0,d-f} \cdot w({\cal O}_{z^2}{\cal O}_{z^{-1}})_{0,f}\nonumber\\
  =&\frac1{2\cdot 3^{d+1}}\prod_{j=0}^d\left(\frac1{2\pi\sqrt{-1}}\oint_{C_j}\frac{dz_j}{z_j^4}\right)z_0(2z_{d-f}-z_{d-f-1}-z_{d-f+1})
\nonumber\\
  &\times\left(\prod_{j=1}^d\tilde{e}(z_{j-1},z_j)\right)\left(\prod_{j=1}^{d-1}\frac1{6z_j(2z_j-z_{j-1}-z_{j+1})}\right)\frac1{z_d}
\label{last} 
\end{align}
 for all $1\leq f\leq d-1$.
Here, $\frac1{2\pi\sqrt{-1}}\oint_{C_j}dz_j$ means taking residues at $z_{0}=0\;(j=0),z_{j}=0,\;\frac{z_{j-1}+z_{j+1}}{2}\;
(j=1,\cdots,d-1),\;z_{d}=0\;(j=d)$.  
\end{lemm}
{\bf Proof.}
Note that $\frac{1}{2}w({\cal O}_{z^1}{\cal O}_{z^0})_{0,d-f}$ is given by, 
\begin{eqnarray*}
&&\frac{1}{2\cdot 3^{d-f+1}}\prod_{j=0}^{d-f}\left(\frac1{2\pi\sqrt{-1}}\oint_{C_j}\frac{dz_j}{z_j^4}\right)z_0\left(\prod_{j=1}^{d-f}
\tilde{e}(z_{j-1},z_j)\right)\\
&&\times\left(\prod_{j=1}^{d-f-1}\frac1{6z_j(2z_j-z_{j-1}-z_{j+1})}\right).
\end{eqnarray*}
Since the integrand in (\ref{last}) is holomorphic at $z_{d-f}=(z_{d-f-1}+z_{d-f+1})/2$, the assertion of the lemma follows from integration of $z_j$'s in ascending order of the subscript $i$.
$\Box$

\vspace{0.5cm}
{\bf Proof of the Main Theorem \ref{theorem:main}.}\\
Let
\begin{equation*}
 \int L_1(e^x)dx=x+\sum_{i=1}^\infty \frac{1}{2}w({\cal O}_{z^1}{\cal O}_{z^0})_{0,d}e^{dx}
\end{equation*}
be a primitive function of $L_1(e^x)$.
Assertion of the theorem \ref{theorem:main} is equivalent to the following equality:
\begin{equation}\label{eq:L0L1}
 f_1(e^x)=f_0(e^x)\int L_1(e^x)dx=L_0(e^x)\int L_1(e^x)dx,
\end{equation}
where we used Lemma \ref{lemm:hol sol}.
Expanding RHS of (\ref{eq:L0L1}), we obtain,
\begin{align*}
 &x\cdot L_0(e^x)+\sum_{d=1}^\infty\left(\frac12 w({\cal O}_{z^1}{\cal O}_{z^0})_{0,d} + \sum_{f=1}^{d-1}(\frac{f}{4} w({\cal O}_{z^1}{\cal O}_{z^0})_{0,d-f} \cdot w({\cal O}_{z^2}{\cal O}_{z^{-1}})_{0,f})\right)e^{dx}.
\end{align*}

We can easily derive the equality:$\sum_{f=1}^{d-1}f(2z_{d-f}-z_{d-f-1}-z_{d-f+1})=d(z_1-z_0)+z_0-z_d$.
Hence application of Lemma \ref{lemm:last} and Lemma \ref{lemm:hol sol} to the RHS of (\ref{eq:L0L1}) results in,
\begin{equation}
 x\cdot f_0(e^x)+\sum_{d=1}^\infty R_d e^{dx},
\end{equation}
where
\begin{align*}
 R_d:=&\frac1{2\cdot 3^{d+1}}\prod_{j=0}^d\left(\frac1{2\pi\sqrt{-1}}\oint_{C_j}\frac{dz_j}{z_j^4}\right)z_0(d(z_1-z_0)+z_0)\\
 &\times\left(\prod_{j=1}^d\tilde{e}(z_{j-1},z_j)\right)\left(\prod_{j=1}^{d-1}\frac1{6z_j(2z_j-z_{j-1}-z_{j+1})}\right)\frac1{z_d}.
\end{align*}

By combining Lemma \ref{lemm:hol sol} and \ref{lemm:second}, we can derive
\begin{equation*}
 R_d=\frac{2^{3d}(6d-1)!!}{(d!)^3}\left(\sum_{j=1}^{3d}\frac{6}{2j-1}-\sum_{j=1}^d\frac3j\right).
\end{equation*}
Therefore, the RHS of  (\ref{eq:L0L1}) becomes,
\begin{equation}
 x\cdot f_0(e^x)+\sum_{d=1}^\infty (\sum_{j=1}^{3d}\frac{6}{2j-1}-\sum_{j=1}^{d}\frac{3}{j})\frac{2^{3d}(6d-1)!!}{(d!)^3}e^{dx}.
\end{equation}
The formula (\ref{eq:solution1}) tells us that it is nothing but $f_{1}(e^{x})$.  $\Box$
\\
\\
\noindent
\begin{rema}
\end{rema}
In the proof of main theorem of Part II, we used the residue argument.
In this remark, we demonstrate how to translate it to Chow ring argument by quoting Lemma \ref{lemm:hol sol}.

We have to compute
\begin{equation}\label{eq:chow comp}
 \int_{\widetilde{Mp}_{0,2}(\bP(1,1,1,3),d)}\frac{H_0^2}{H_d}E_d,
\end{equation}
where.
\begin{align*}
 E_d:&=\frac{\prod_{i=1}^de^6(H_{i-1},H_i)}{\prod_{i=1}^{d-1}6H_i}\\
 &=\frac1{3^{d+1}}\cdot\frac{\prod_{i=1}^d\tilde{e}(H_{i-1},H_i)}{\prod_{i=1}^{d-1}6H_i}\cdot \left(3^{d+1}\prod_{i=0}^{d-1}(2H_i+H_{i+1})(H_i+2H_{i+1})\right).
\end{align*}
Since $E_d$ has a factor $3^{d+1}\prod_{i=0}^{d-1}(2H_i+H_{i-1})(H_i+2H_{i+1})$, we can compute (\ref{eq:chow comp}) as
\begin{align*}
 \frac{H_0^2}{H_d}\frac1{3^{d+1}}\cdot\frac{\prod_{i=1}^d\tilde{e}(H_{i-1},H_i)}{\prod_{i=1}^{d-1}6H_i}
\end{align*}
in $\C[H_0,\dots,H_d]/(H_0^4,H_1^4(-H_0+2H_1-H_2),\dots,H_{d-1}^4(-H_{d-2}+2H_{d-1}-H_d),H_d^4)$.

We can prove that
\begin{equation*}
 H_0^3H_1^4H_2^4\cdots H_{j-1}^4H_j^5=\frac{j}{j+1}H_0^3H_1^4H_2^4\cdots H_j^4H_{j+1}
\end{equation*}
for all $j=0,1,2,\dots,d-1$ by mathematical induction.
Then, we obtain
\begin{align*}
 &\frac{H_0^3H_1^4H_2^4\cdots H_{j-2}^4H_{j-1}^3\tilde{e}(H_{j-1},H_j)}{H_j}\\
 =&2^4\cdot 3^2 H_0^3H_1^4H_2^4\cdots H_{j-2}^4H_{j-1}^4(5H_{j-1}+H_j)(3H_{j-1}+3H_j)(H_{j-1}+H_j)\\
 =&2^4\cdot 3^2 H_0^3H_1^4H_2^4\cdots H_{j-2}^4H_{j-1}^4(\frac{5(j-1)}{j}+1)(\frac{3(j-1)}{j}+3)(\frac{j-1}{j}+5) H_j^3\\
 =&2^4\cdot 3^2\frac{(6j-1)!!}{j^3\cdot(6j-7)!!} H_0^3H_1^4H_2^4\cdots H_{j-2}^4H_{j-1}^4H_j^3.
\end{align*}
This gives us the above successive argument.


\newpage

\begin{thebibliography}{9}
  \bibitem{CdGP} P.Candelas, X.de la Ossa, P.Green and L.Parkes, \textit{A pair of Calabi-Yau manifolds as an exactly soluble superconformal field theory}, Nuclear Physics {\bf B395} (1991).
\bibitem{CK1} I. Ciocan-Fontanine, B. Kim. \textit{Wall-crossing in genus zero quasimap theory and mirror maps},  Algebr. Geom. 1 (2014), no. 4, 400--448. 
\bibitem{CCK} D. Cheong, I. Ciocan-Fontanine, B. Kim, \textit{Orbifold quasimap theory}, Math. Ann. 363 (2015), no. 3--4, 777--816. 
\bibitem{CK} D. A. Cox, S.Kats, \textit{Mirror Symmetry and Algebraic Geometry}, Mathematical Surveys and Monographs Vol.68, American Mathematical Society,1999.
 \bibitem{Cox} D. A. Cox, \textit{The homogeneous coordinate ring of a toric variety, J. Algebraic Geometry}, {\bf 4} (1995), 17-50, alg-geom/9210008.
  \bibitem{Ful} W. Fulton, \textit{Introduction to Toric Varieties}, Princeton University Press, Princeton, 1993.
  \bibitem{Givental} A. Givental, \textit{A mirror theorem for toric complete intersections}, Topological field theory, primitive forms and related topics (Kyoto, 1996), 141-175, Progr. Math., 160, Birkh\"{a}user Boston, Boston, MA, 1998.
  \bibitem{Jin1} M. Jinzenji, \textit{Mirror Map as Generating Function of Intersection Numbers: Toric Manifolds with Two K\"{a}hler Forms}, arXiv:1006.0607, Comm. Math. Phys. 323(2013), no. 2, 747-811.
 \bibitem{LLY} B. Lian, K. Liu and S. T. Yau, \textit{Mirror Principle III}, Asian J. Math. {\bf 3} (1999),no.4, 771-800
  \bibitem{LY} B. Lian, S. T. Yau, \textit{Arithmetic Properties of Mirror Map and Quantum Coupling}, Comm. Math. Phys. 176(1996), no. 1, 163-191.
  \bibitem{Sch} S. Scholtes, \textit{Introduction to Piecewise Differentiable Equations}, Springer, 2012.
  \bibitem{S} H. Saito, \textit{Chow Rings of $\widetilde{Mp}_{0,2}(N,d)$ and $\overline{M}_{0,2}(\bP^{N-1},d)$ and Gromov-Witten Invariants of Projective Hypersurfaces of Degree 1 and 2}, Internat. J. Math. 28(2017), no. 12, 1750090.
\end{thebibliography}
\end{document}